\documentclass{article}
\usepackage[utf8]{inputenc}
\usepackage[hidelinks]{hyperref}
\usepackage{graphicx} % Required for inserting images
\usepackage[margin=1.5in]{geometry}
\usepackage{amsthm}
\usepackage{times}
\usepackage{newtxmath}
\usepackage{fontawesome5}
\usepackage{xcolor}
\usepackage{array}
\usepackage[natbib=true,style=authoryear,maxcitenames=2,maxbibnames=6,]{biblatex}
\usepackage[createShortEnv]{proof-at-the-end}
\usepackage{authblk}

\newcommand{\M}{\mathcal M}
\newcommand{\eps}{\varepsilon}

\newcommand{\Sset}{\mathbb S}
\newcommand{\Spairs}{\Sset_\text{pairs}}
\newcommand{\D}{\mathbb D}
\newcommand{\G}{\mathcal G}
\newcommand{\R}{\mathbb R}
\newcommand{\I}{\mathcal I}

\newcommand{\range}{\mathrm{range}}

\renewcommand{\Pr}{\mathbf{P}}

\newtheorem{theorem}{Theorem}

\theoremstyle{definition}
\newtheorem{definition}[theorem]{Definition}

\theoremstyle{remark}
\newtheorem{example}{Example}
\newtheorem{remark}{Remark}

\title{Interpreting Network Differential Privacy}
\author{Jonathan Hehir}
\author{Xiaoyue Niu}
\author{Aleksandra Slavkovi{\'c}}
\affil{Department of Statistics, Pennsylvania State University, University Park, PA}
\date{\today}

\begin{document}

\maketitle

\begin{abstract}
    How do we interpret the differential privacy (DP) guarantee for network data? We take a deep dive into a popular form of network DP ($\eps$--edge DP) to find that many of its common interpretations are flawed. Drawing on prior work for privacy with correlated data, we interpret DP through the lens of adversarial hypothesis testing and demonstrate a gap between the pairs of hypotheses actually protected under DP (tests of complete networks) and the sorts of hypotheses implied to be protected by common claims (tests of individual edges). We demonstrate some conditions under which this gap can be bridged, while leaving some questions open. While some discussion is specific to edge DP, we offer selected results in terms of abstract DP definitions and provide discussion of the implications for other forms of network DP.
\end{abstract}

\section{Introduction}

In recent years, differential privacy (DP) has emerged as the preeminent formal definition of data privacy, with an explosion of both theoretical developments and real-world, large-scale implementations. Through simple and general probabilistic constraints imposed on a data release mechanism, DP offers a strong and quantifiable privacy guarantee that is adaptable to data of various formats and agnostic to the underlying distribution of potentially sensitive data. Broadly speaking, DP mechanisms are designed to allow population-level inference from a given database while hiding the potential contribution of any individual to that database.

The simplicity, generality, and adaptability of differential privacy are highly desirable qualities, but these same qualities belie considerable complexity in the application and interpretation of the privacy definition, particularly in data with complex, interdependent structures. In this paper, we focus on applications of DP to networks, which are typified by dependence between records and dyadic relationships that leave unclear which data ``belong'' to a given individual.

Network DP comes in several variants, including the common edge DP and node DP \citep{hay2009accurate} as well as a number of less common definitions that are tailored to the specific privacy risk, network structure, and inferential task relevant to a given problem. Proper interpretation of these DP variants and their implications for individuals' privacy are important both for our understanding as practitioners as well as for effective communication of the privacy guarantee to constituents whose data are subject to DP.

In the literature on network privacy, we find a concerning pattern of inaccurate or misleading interpretations of DP that suggest a stronger notion of privacy than actually afforded. For example, a common claim is that edge DP (Definition~\ref{def:edge-dp}) bounds the ability of an adversary to infer the presence of any single edge in the network. In reality, such a claim only holds under extremely restrictive assumptions on the underlying distribution of the network or the knowledge of the adversary. These assumptions are generally unrealistic but frequently left unstated. (See Section~\ref{sec:interpretations}.) A major goal of our work is to determine reasonable conditions under which similar claims can be made.

\textbf{Related Work.} \citet[Appendix~A]{tschantz2020sok} provides a detailed history of varied perspectives on the interpretation of DP, noting two dominant trends. One approach is to quantify what can be learned about individuals from a DP release. These \textit{inferential-privacy} interpretations follow from DP in the presence of additional assumptions on either the data or the adversary, which has often led to erroneous statements that DP in fact requires such assumptions. For example, \citet{liu2016dependence} states that DP requires independent-tuple databases, and \citet{cuff2016differential} states that DP requires the adversary to know all but one record in the database. It has long been shown that general limits on inferential privacy are impossible and argued that DP does not (in general) imply these results \citep{dwork2010difficulties,ghosh2016inferential}. A number of works have been dedicated to demonstrating that DP does not protect against various forms of sensitive inference \citep[e.g.,][]{kifer2009attacks,cormode2011personal,kifer2014pufferfish} or to deriving bounds on adversarial inference under various constraints \citep[e.g.,][]{kifer2014pufferfish,chen2014correlated,yang2015bayesian,liu2016dependence,cuff2016differential,ghosh2016inferential}.

A second way to interpret DP is as a guarantee that individual contributions to a database are hidden. This is what \citet{tschantz2020sok} refers to as a \textit{causal} view of privacy: the effect of one individual contributing their data point to a database should have minimal effect on the DP release. While ``a bad disclosure can still occur, [DP] assures the individual that it will not be the presence of her data that causes it, nor could the disclosure be avoided through any action or inaction on the part of the user'' \citep{dwork2006differential}. In the network setting, however, what it means for an individual to contribute their data point(s) is often unclear or misaligned with the privacy definition in use. In practice, inferential-privacy interpretations appear to be common in network DP (see examples in Section~\ref{sec:interpretations}).

The field of DP is host to a large number of privacy definitions: \citet{desfontaines2020sok} provides a general overview, while \citet{mueller2022sok} focuses specifically on the variants applied to networks. In this work, we focus primarily on the commonly used definition of edge DP from \citet{hay2009accurate}. A recent survey of applications of DP to network data is given in \citet{jiang2021applications}.

This work relies heavily on the Pufferfish privacy definition of \citet{kifer2014pufferfish} to assess the inferential-privacy guarantees that are implied by edge DP. Pufferfish is a generalization of DP that accommodates the protection of arbitrary secrets, requiring as a trade-off the specification of possible data-generating distributions. We provide additional background on Pufferfish and its relationship with DP in Section~\ref{sec:pufferfish}.

\textbf{Outline.} We review the  classical (non-network) setting of differential privacy in Section~\ref{sec:dp}, providing discussion on both inferential and causal interpretations.
Section~\ref{sec:edge-dp} offers a review of edge DP and its interpretations, where we find causal interpretations are difficult, and common inferential interpretations are flawed.
We turn to Pufferfish in Section~\ref{sec:pufferfish} as a framework to formalize both the forms of inference we wish to protect and the corresponding assumptions that are required to protect it. We demonstrate conditions under which edge DP implies meaningful, quantifiable inferential-privacy guarantees.
We conclude with a discussion of other forms of network DP, the general implications of our work, and remaining open questions in Section~\ref{sec:conclusion}.

\textbf{Notation.} A graph is a collection of vertices and edges, $G = (V, E)$, where $E \subseteq \I_V \subseteq V \times V$, for $\I_V$ the set of all possible edges. Although our work generalizes to directed or undirected graphs, our notation assumes undirected edges and uses the notation $\{i, j\} = \{j, i\}$ to denote the unordered edge between nodes $i$ and $j$. We typically assume the vertex set $V = [n] = \{1, \dots, n\}$ and so occasionally use $G$ as a shorthand for the edge set. Related shorthands include $G \cup \{i, j\} = (V, E \cup \{i, j\})$ and $G \setminus \{i, j\} = (V, E \setminus \{i, j\})$. We sometimes use lowercase notation for graphs to distinguish a fixed graph $g$ from a random graph $G$. The operator $\triangle$ denotes symmetric difference, i.e., $A \triangle B = (A \setminus B) \cup (B \setminus A)$. $\lnot$ is the logical \textit{not} operator. %$1_\text{predicate}$ is an indicator variable that takes the value 1 if \textit{predicate} is true and 0 otherwise.

\section{Review of Differential Privacy}
\label{sec:dp}

In its most general form, differential privacy \citep{dwork2006calibrating} is a constraint applied to a randomized algorithm (or \textit{privacy mechanism}) $\M : \D \to \range(\M)$ acting on a \textit{database} from some space of databases $\D$ equipped with a symmetric \textit{neighboring} relationship $\sim$. In classical settings concerning tabular data, $\D$ is usually chosen to be the class of databases consisting of lists of tuples (rows) from some space of tuples $\mathbb T$, and the neighboring relationship $\sim$ holds if two databases differ on one (and only one) tuple. (We will formalize this notion below.) Specifically, DP requires that the distribution of $\M$ applied to any neighboring databases must be close in the following sense.

\begin{definition}[Differential Privacy, \citealp{dwork2006calibrating}]
\label{def:dp}
Suppose $\eps > 0$, $\D$ is a class of databases, and $\sim$ is a symmetric relation on $\D$. Let $\M : \D \to \range(\M)$ be a randomized algorithm. We say $\M$ satisfies \textbf{$\eps$--DP$(\D, \sim)$} if for all $D_0, D_1 \in \D$ satisfying $D_0 \sim D_1$ and all $R \subseteq \range(\M)$,
\begin{equation}
\label{eq:general-dp}
\Pr(\M(D_0) \in R) \leq e^\eps \, \Pr(\M(D_1) \in R) .
\end{equation}
\end{definition}

Definition~\ref{def:dp} presents DP in a general form that captures essentially any form of \textit{pure differential privacy}%
\footnote{For simplicity, we restrict our attention to pure DP to avoid the interpretative complications of approximate DP forms such as $(\eps, \delta)$-DP. See discussion in Section~\ref{sec:conclusion}.}
through the appropriate selection of $\D$ and $\sim$. We emphasize that $D_0$ and $D_1$ are treated as deterministic in this definition and that the probability measure $\Pr$ in Eq.~\ref{eq:general-dp} captures only the randomness of the release mechanism $\M$. The constraint is applied to every pair of neighboring databases $D_0 \sim D_1$, without respect to the probability (or improbability) of their occurrence. Thus by definition, a DP guarantee is robust to arbitrary underlying data-generating distributions.

Recall that in classical tabular data, $\D$ is chosen to be the class of databases consisting of lists of tuples from $\mathbb T$. The neighboring relationship $\sim$ comes in two main varieties:
\begin{itemize}
    \item \textit{Bounded neighbors:}
        $D_0 \sim D_1$ iff $|D_0| = |D_1| = n$ and $D_0, D_1$ agree on all but one tuple.%
        \footnote{Sometimes referred to as \textit{$\eps$-indistinguishability} in early literature.}
    \item \textit{Unbounded neighbors:}
        $D_0 \sim D_1$ iff $|D_0| = |D_1| + 1$ and $D_1 \subset D_0$ (or vice-versa).
\end{itemize}
Under a bounded-neighbors definition, the size of the database (i.e., number of rows, $n$) is fixed and not sensitive, and the privacy guarantee essentially states that if an individual record is changed, the distribution of $\M$ should not substantially change. By contrast, with unbounded neighbors, the size of the database is sensitive, and the output of $\M$ should not substantially change if any one record is excluded from the database.

Since each record in the classical database is typically assumed to belong to an individual, the concept of neighboring databases captures pairs of hypothetical databases that would arise from a given individual contributing truthful information vs. contributing fabricated information (bounded neighbors) or not contributing any information at all (unbounded neighbors). In this setting, DP may be given a more approachable interpretation: \textit{The output of the DP release $\M(D)$ is not substantially changed (as quantified by $\eps$) whether or not a given individual's data is (faithfully) included in $D$.} This is what \citet{tschantz2020sok} refers to a \textit{causal} interpretation of privacy. In the following section, we will argue that causal interpretations generally do not apply as gracefully to network DP variants.

Interpretations of DP can be given statistical flavor by stating results in terms of hypothesis testing. Because the DP constraint requires that the released data $\M(D)$ be approximately indistinguishable between neighboring databases, an adversary observing $\M(D)$ should not be able to reliably detect which specific database $D$ was processed by the algorithm $\M$. We formalize this in Theorem~\ref{thm:dp-distingishing-neighbors}, which restates a well-known result from \citet{kairouz2015composition}.

\begin{theoremE}[\citealp{kairouz2015composition}]
    \label{thm:dp-distingishing-neighbors}
    Let $\M$ be a randomized algorithm, and let $D_0 \sim D_1$ in $\D$. Then $\M$ is $\eps$-DP$(\D, \sim)$ if and only if any level-$\alpha$ hypothesis test for $H_0: D = D_0$ vs. $H_1: D = D_1$ that is a function of $\M$ and $\M(D)$ has power $\beta \leq e^\eps \alpha$.
\end{theoremE}

\begin{proofE}
    This follows from \citet[Theorem~2.1]{kairouz2015composition}, but for completeness we provide a proof of our restatement.
    
    ($\implies$) Assume $\M$ is $\eps$-DP. By the post-processing invariance of DP \citep{dwork2006calibrating}, $\mathcal T \circ \M$ satisfies $\eps$-DP($\D, \sim$), and so for all $D_0 \sim D_1$ from $\D$:
    \begin{align*}
        \text{$\mathcal T$ is level-$\alpha$} &\implies \Pr( (\mathcal T \circ \M)(D_0) = 1) \leq \alpha \\
            &\implies \beta = \Pr( (\mathcal T \circ \M)(D_1) = 1 ) \leq e^\eps \alpha .
    \end{align*}

    ($\impliedby$) On the other hand, assume the hypothesis-testing condition holds, and let $E \subseteq \range(\M)$ be any event. Choosing $\mathcal T(\omega) = 1_{\omega \in E}$, $\alpha = \Pr(\M(D_1) \in E)$, we have that $\mathcal T$ is level-$\alpha$ for the test $H_0: D = D_1$ vs. $H_1: D = D_0$. Since the power of this test is bounded by $e^\eps \alpha$, we have $\Pr(\M(D_0) \in E) \leq e^\eps \Pr(\M(D_1) \in E)$. Since this holds for any $D_0 \sim D_1$ from $\D$ and $E \subseteq \range(\M)$, the condition of $\eps$-DP($\D, \sim$) is satisfied.
\end{proofE}
The hypothesis test depicted in Theorem~\ref{thm:dp-distingishing-neighbors} captures the essential statistical promise of DP, but this test may feel a bit unnatural, treating the entirety of the database as the parameter of interest. Not only is this unusual statistical practice, it may seem to capture only a very specific type of privacy violation: the ability to infer a sensitive database in its entirety. Couldn't inference about a mere portion of the database---one individual's data, let's say---constitute a privacy violation? The answer to this turns out to be somewhat complicated, prone to both confusion and controversy.% \citep[see, e.g.,][]{mcsherry2016github,mcsherry2016githublunch}.

As a concrete example, consider again the classical tabular setting under the bounded-neighbors definition of DP. Suppose that a database consists of a collection of records $D = (X_1, \dots, X_n)$ corresponding to $n$ individuals, and suppose the adversary knows the record values for all but the $i$-th individual. In this case, we can indeed say that the adversary observing $\M(D)$ cannot test $H_0: X_i = t_0$ vs. $H_1: X_i = t_1$ for any tuples $t_0, t_1$ at level $\alpha$ with power greater than $e^\eps \alpha$. After all, conditioned on the remaining $n-1$ records, the hypothesis tests for the $i$-th record are equivalent to hypothesis tests for the full database. This hypothetical almost-all-knowing adversary is often referred to as a \textit{strong adversary}. Curiously, if we assume instead a weaker adversary, such a statement no longer holds \citep{kifer2009attacks,kifer2014pufferfish}. To make a comparable statement about adversarial hypothesis testing at the individual-record level for weaker adversaries requires further qualification. Perhaps the most well known such statement is the one provided by \citet[Theorem~2.4]{wasserman2010statistical}, which assumes the $X_i$ tuples are independent and identically distributed. We restate this in slightly more general terms here, removing the requirement of identical distribution.

\begin{theoremE}[\citealp{wasserman2010statistical}]
    \label{thm:wasserman}
    Let $D = (X_1, \dots, X_n)$ be independent random tuples in $\mathbb T_1 \times \dots \times \mathbb T_n$. Let $\M$ be $\eps$-DP$(\mathbb T_1 \times \dots \times \mathbb T_n, \sim)$ for $\sim$ capturing bounded-neighbors as defined previously. For all $i \in [n]$ and $t_0, t_1 \in \mathbb T_i$, any level-$\alpha$ hypothesis test for $H_0: X_i = t_0$ vs. $H_1: X_i = t_1$ that is a function of $\M$, $\M(D)$, and the distribution functions of $X_1, \dots, X_n$ has power $\beta \leq e^\eps \alpha$.
\end{theoremE}

\begin{proofE}
    This is a minor restatement of Theorem~2.4 of \citet{wasserman2010statistical}. Note that while the theorem statement in \citet{wasserman2010statistical} does not explicitly mention independence, the $X_i$ are previously defined to be both independent and identically distributed. In fact, the proof they provide holds under the slightly more general assumption of independent (but not necessarily identically distributed) $X_i$.
\end{proofE}
Here we arrive at what may look like a contradiction: the definition of DP (Definition~\ref{def:dp}) requires neither assumptions about the distribution of sensitive data, nor assumptions about the potential knowledge of the adversary. These are two of DP's key strengths. On the other hand, interpretations of DP that wrestle with the ability to infer an individual's data (e.g., Theorem~\ref{thm:wasserman}) require precisely these sorts of assumptions. Returning to the language of \citet{tschantz2020sok}, while the causal interpretation of DP always holds, these \textit{inferential-privacy} interpretations of DP are known to require additional assumptions \citep{dwork2010difficulties}. 

In a series of blog posts \citep{mcsherry2016github,mcsherry2016githublunch}, Frank McSherry approaches this nuance by drawing a distinction between what he calls \textit{your secrets} and \textit{secrets about you}. \textit{Your secrets} are the secrets that can only be learned if you share information about them. The causal interpretation of DP implies the protection of your secrets. \textit{Secrets about you}, on the other hand, could be any sensitive information concerning you. Inferential privacy aims to protect both types of secrets. While DP does protect \textit{your secrets}, it offers no general guarantee for all \textit{secrets about you}. This eloquent distinction is helpful in the classical tabular setting, but attempting to grapple with network DP through this framework may prove challenging, as the definitions of \textit{your secrets} and \textit{secrets about you} blur. A relationship between you and me is a secret about either of us, but is it your secret, my secret, our secret, or some other option? And given network patterns like reciprocity and transitivity, what does one network secret reveal about another?

\section{The Case of Edge Differential Privacy}
\label{sec:edge-dp}

A clear definition of neighboring databases is imperative in any application of DP. In classical settings, this often reduces to the bounded vs. unbounded neighbors distinction drawn in Section~\ref{sec:dp}. In network settings, however, there are many ways to define neighboring networks. Perhaps the most popular of these is \textit{edge differential privacy}, which, loosely speaking, defines two networks as neighboring if they differ on a single edge. We explore here some of the subtleties of edge DP, its interpretations and common misinterpretations, and the implications of edge DP on adversarial inference.

Consider a class of binary, undirected networks. In general, we will denote a network as a graph $G = (V, E)$, where we assume the vertex set $V = [n]$ and edge set $E \subseteq \I_V$, where $\I_V$ is the set of valid edges for the graph. The discussion to follow largely generalizes to graphs with directed or undirected edges, weighted or unweighted edges, and with or without self-loops. Let $\G_n$ denote the set of undirected, binary graphs with $n$ vertices.

In the original specification of edge DP \citep{hay2009accurate}, neighboring networks differ on a single edge but do not need to have the same number of nodes---an unbounded flavor, if you will. Specifically, Hay's edge DP allows for any single edge to be changed \textit{or for the addition or removal of an isolated node}. Many popular adaptations of edge and node DP, however, take a fixed-size approach, requiring the DP constraint only hold over networks with the same number of nodes \citep[e.g.,][]{karwa2011private,kasiviswanathan2013analyzing,karwa2016inference}. %\todo{More citations?}
This distinction is particularly important if the presence or absence of a node in the network is sensitive. The bounded flavors of these variants are often accompanied by an explicit assumption that certain information (e.g., node count, node identities) is assumed to be publicly known \citep[e.g.,][]{kasiviswanathan2013analyzing,karwa2017sharing,hehir2021consistent}. 
%\todo{More/better example citations.}
We provide a formal definition in bounded flavor below:

\begin{definition}[Edge DP]
    \label{def:edge-dp}
    Let $G_1 = (V^{(1)}, E^{(1)}), G_2 = (V^{(2)}, E^{(2)}) \in \G_n$.
    We write $D_1 \overset{e}{\sim} D_2$ iff $V^{(1)} = V^{(2)}$ and $|E^{(1)} \triangle E^{(2)}| = 1$, where $\triangle$ denotes symmetric difference (i.e., only one edge differs between $G_1$ and $G_2$).
    The randomized algorithm $\M : \G_n \to \range(\M)$ satisfies \textbf{$\eps$--edge DP} if it satisfies $\eps$-DP$(\G_n, \overset{e}{\sim})$.
\end{definition}

\subsection{Interpreting (and Misinterpreting) Edge DP}
\label{sec:interpretations}

Many works that propose edge-DP methods stop short of providing clear interpretations of the privacy guarantee. We speculate that this is due to the sheer difficulty of providing an interpretation that is clearly intelligible and relatable to actual privacy concerns. The most direct translation of the edge DP definition (Definition~\ref{def:edge-dp}) into plain English reads like this: if an algorithm is edge-DP, then changing a single edge in the network should not substantially change the distribution of possible outputs. This, however, says very little---at least directly---about privacy.

Alternatively, we can say that given two hypothetical networks, one in which a single edge is included and another otherwise identical network in which that same edge is not included, an adversary observing an edge-DP output produced from either of these two networks should reach approximately the same conclusions, not learning considerably more from one than the other. While this interpretation incorporates an adversary, it is unclear to what situations these two hypothetical networks correspond.

Fundamentally, these causal-flavored interpretations feel unnatural because there often lacks a natural setting corresponding to the notion of neighboring databases in edge DP. Recall that in classical DP, neighboring databases capture the distinction between an individual contributing true data vs. fabricated data (bounded neighbors) or contributing true data vs. witholding data (unbounded neighbors). These are the very scenarios that immediately occur to individuals concerned about their privacy in data processing. But in what hypothetical situations do two networks arise that differ on a single edge, and how do these situations correspond to privacy? To whom is a given network edge sensitive, and who makes the decision to report or not report a link between two individuals?

Avoiding the awkwardness of neighboring databases in edge DP has led to inferential-flavored interpretations that are more pleasing to the ear but lack crucial detail. Consider, for example, the interpretation given with the original definition of edge DP \citep[][emphasis added]{hay2009accurate}:
\begin{quote}
    Clearly, any [differentially private] observed output cannot reveal much about [a] particular record if that output is (almost) as likely to occur even when the record is excluded from the database. [....] %To adapt differential privacy to graphs, we must choose a definition for neighboring graphs and understand the privacy semantics of that choice. [....] 
    \textit{An edge--differentially private algorithm protects individual edges from being disclosed.}
\end{quote}
Variations on the short sentence in italics are abundant in the literature on network DP, suggesting---often explicitly---that edge DP guarantees to limit the ability of an adversary to detect or infer a single edge in the network. Examples span the range from old to new and include a number of influential works \citep[e.g.,][]{task2012guide,blocki2013differentially,task2014should,xiao2014differentially,karwa2016inference,karwa2017sharing,jiang2021applications}. Sadly, while these concise interpretations sound more satisfying than those given above, this powerful simplicity proves too good to be true, as demonstrated next.

\subsection{Detecting Edges under Edge Privacy}

We illustrate below a case where an $\varepsilon$--edge DP mechanism fails to protect against inferring the presence or absence of a specific edge in a network. We reiterate that this is not a failing of the privacy guarantee---edge DP does not guarantee protection against arbitrary edge-level inference---but rather a counterexample to a common misinterpretation of edge DP's guarantee.

To start, we need a mechanism that satisfies $\varepsilon$--edge DP. For simplicity, we employ a mechanism that releases an estimate of the number of edges in a graph using the well-known Laplace mechanism.

\begin{definition}[Laplace Mechanism, \citealp{dwork2006calibrating}]
    Let $\D$ be a class of databases, $\sim$ be a symmetric relation on $\D$, and $f : \D \to \R$ be a function. Define the sensitivity $\Delta = \sup_{D_0 \sim D_1} |f(D_0) - f(D_1)|$ over all neighboring databases $D_0, D_1 \in \D$. The $\varepsilon$--DP$(\D, \sim)$ Laplace mechanism $\M : \D \to \R$ is given by:
    \[
    \M_\text{Laplace}(D) = f(D) + \mathrm{Laplace}(\Delta/\varepsilon) .
    \]
\end{definition}

Adapting the Laplace mechanism to the problem of counting edges in $\varepsilon$--edge DP is straightforward. By Definition~\ref{def:edge-dp}, $|E^{(0)} \triangle E^{(1)}| = 1$ for any $G_0 \overset{e}{\sim} G_1 \in \G_n$. This means the sensitivity of the edge count is exactly one, i.e., $\Delta = \sup_{G_0 \overset{e}{\sim} G_1} \left| \, | E^{(0)} | - | E^{(1)} | \, \right| = 1$. As a result, the function $\M_\text{edges}: \G_n \to \R$ that releases a noisy count of edges from the graph $G = (V, E) \in \G_n$,
\[
\M_\text{edges}(G) = |E| + \mathrm{Laplace}(1 / \varepsilon) ,
\]
satisfies $\varepsilon$--edge DP.

Having established a privacy mechanism, we now provide an example in which this mechanism fails to prevent inference of the presence or absence of a specific edge.

\begin{example}
\label{ex:insects}
Suppose $G = ([n], E)$ denotes a hive of $n > 2$ insects, where insects $1$ and $2$ are the two queens of the hive. An undirected edge $\{i, j\} \in E$ denotes that insects $i$ and $j$ have a cooperative working relationship. When the two queens are cooperating (i.e., $\{1, 2\} \in E$), all other pairs of insects cooperate independently with probability $a \in [0, 1]$. When the two queens are not cooperating, all other pairs cooperate independently with probability $b \in [0, a)$. Let $\tilde d = {n \choose 2}^{-1} \M_\text{edges}(G)$. Then $\tilde d$ satisfies $\eps$--edge DP, while $1_{\tilde d > (a+b)/2}$ converges in probability to $1_{\{1, 2\} \in E}$ as $n \to \infty$.
\end{example}

Example~\ref{ex:insects} releases a noisy estimate of the network density that satisfies $\eps$--edge DP. However, when $n$ is large, this density estimate can be used to infer the existence of a cooperative relationship between the two queens with high probability, since the noisy density of the graph tends to $a$ when the two queens cooperate vs. $b$ if they do not cooperate. In short, it is virtually impossible to hide the existence of an edge between the queens simply by masking the existence of that single edge. Since the relationships of the remaining insects depend on the relationship between the queens, the secret of the queens' relationship is conflated with all the remaining relationships in the hive.

While this is a dramatic simplification of real network dynamics---illustrated on a species not generally subject to reasonable expectations of privacy---the lesson is the same for real-world networks. Networks are characterized by dyadic relationships and forms of dependence that make muddy concepts like \textit{your secrets}. In networks representing common patterns such as reciprocity, transitivity, and homophily, any given relationship (or lack thereof) may reveal information about another potential relationship, and every relationship (or absence thereof) involves not one but two individuals.

In short, the traditional interpretations of DP reviewed in Section~\ref{sec:dp} don't extend gracefully to edge DP. Causal interpretations struggle to explain useful privacy hypotheticals, leading to an abundance of inferential misinterpretations.
While troublesome without qualification, inferential interpretations can be made in the presence of specific assumptions about the adversary's knowledge or the underlying network distribution. In light of the popularity of these statements and the difficulties associated with causal interpretations, we turn our attention now to a brief exploration of the assumptions required to state more general adversarial bounds under edge DP.

\section{Interpretable Privacy through Pufferfish}
\label{sec:pufferfish}

We return to hypothesis testing as our framework for defining the limits of adversarial inference under privacy. Recall from Theorem~\ref{thm:dp-distingishing-neighbors} that DP protects against a very specific form of hypothesis testing: tests between two complete, neighboring databases. Revisiting the insect example from earlier (Example~\ref{ex:insects}), an adversary who knows (or otherwise assumes) every relationship except the queens' cannot tell whether $\tilde d$ was calculated from a network in which the queens' relationship was marked as present or absent. These two possibilities correspond to two hypothetical, complete databases that could have been used to calculate $\tilde d$---but only one of these databases is plausible under the constraints of the network data generating process. This hypothesis test is essentially meaningless to an adversary knowledgeable of the network constraints.

What we'd like instead is a framework for protecting against sets of hypotheses that are meaningful in these sorts of settings. This is exactly what Pufferfish \citep{kifer2014pufferfish} affords: a framework that generalizes DP but protects against flexible pairs of hypotheses. This flexiblity will allow us to explicitly define the secrets we want to protect in a network and the assumptions required to achieve those protections. We define Pufferfish in Definition~\ref{def:pufferfish} below and provide an interpretation in terms of hypothesis tests in Theorem~\ref{thm:pufferfish-hypothesis-tests}.

\begin{definition}[Pufferfish privacy, \citealp{kifer2014pufferfish}]
\label{def:pufferfish}
Given a set of potential secrets $\Sset$, a set of discriminative pairs $\Spairs$, a set of data-generating distributions $\Theta$, and a privacy parameter $\eps > 0$,
a (potentially randomized) algorithm $\M$ satisfies \textbf{$\eps$--Pufferfish$(\Sset, \Spairs , \Theta)$ privacy} if
\begin{enumerate}
    \item for all possible outputs $\omega \in \range(\M)$,
    \item for all pairs $(s_i , s_j ) \in \Spairs$ of potential secrets, and
    \item for all $\Pr_\theta \in \Theta$ for which $\Pr_\theta(s_i) \neq 0$ and $\Pr_\theta(s_j) \neq 0$
\end{enumerate}
the following hold for $D \sim \Pr_\theta$:
\[
\begin{aligned}
\Pr_\theta(\M(D) = \omega \mid s_i) &\leq e^\eps \Pr_\theta(\M(D) = \omega \mid s_j) \\
\Pr_\theta(\M(D) = \omega \mid s_j) &\leq e^\eps \Pr_\theta(\M(D) = \omega \mid s_i) .
\end{aligned}
\]
\end{definition}

\begin{remark}
As in \citet{kifer2014pufferfish}, we present Pufferfish in terms of discrete probability distributions. Continuous distributions may be accommodated through the appropriate use of continuous probability density functions.
\end{remark}

\begin{theoremE}[][end]
    \label{thm:pufferfish-hypothesis-tests}
    Let $\M$ be a randomized algorithm, $\Theta$ a set of data-generating distributions, $\Sset$ a set of potential secrets, and $\Spairs$ a set of discriminative pairs. Then the following two statements are equivalent:
    \begin{itemize}
        \item $\M$ is $\eps$-Pufferfish($\Sset, \Spairs, \Theta$).
        \item For any database $D \sim \Pr_\theta$ for $\Pr_\theta \in \Theta$ and all choices of $\alpha > 0$ and $(s_0, s_1) \in \Spairs$ satisfying $\min \{ \Pr_\theta(s_0), \Pr_\theta(s_1) \} > 0$, any level-$\alpha$ hypothesis test of $H_0: s_0$ vs. $H_1: s_1$ (or $H_0: s_1$ vs. $H_1: s_0$) that is a function of $\M$, $\M(D)$, and $\Pr_\theta$ has power $\beta \leq e^\eps \alpha$.
    \end{itemize}
\end{theoremE}

\begin{proofE}
    ($\implies$) Assume the Pufferfish condition is satisfied, and fix $\Pr_\theta \in \Theta$ such that $D \sim \Pr_\theta$. Without loss of generality, consider the test for $H_0: s_0$ vs. $H_1: s_1$. By the post-processing invariance of Pufferfish \citep{kifer2014pufferfish}, $\mathcal T \circ \M$ satisfies $\eps$-Pufferfish($\Sset, \Spairs, \D$), and so for any $(s_0, s_1) \in \Spairs$ satisfying $\min \{ \Pr_\theta(s_0), \Pr_\theta(s_1) \} > 0$:
    \begin{align*}
        \text{$\mathcal T$ is level-$\alpha$} &\implies \Pr_\theta( (\mathcal T \circ \M)(D) = 1 \mid s_0 ) \leq \alpha \\
            &\implies \beta = \Pr_\theta( (\mathcal T \circ \M)(D) = 1 \mid s_1 ) \leq e^\eps \alpha .
    \end{align*}

    ($\impliedby$) Assume the Pufferfish condition is \textit{not} satisfied. Then without loss of generality in the order of $s_0$ and $s_1$, there exists a distribution $\Pr_\theta \in \Theta$, a discriminative pair $(s_0, s_1) \in \Spairs$, and an event $R \subseteq \range(\M)$ such that:
    \[
    \Pr_\theta( \M(D) \in R \mid s_0) > e^\eps \Pr_\theta( \M(D) \in R \mid s_1) .
    \]
    Let $\alpha = \Pr_\theta( \M(D) \in R \mid s_0)$, and let $\mathcal T(\omega) = 1_{\omega \in R}$.
    Then $\alpha > 0$ and $\mathcal T$ is a level-$\alpha$ test, since
    \[
    \Pr_\theta( \mathcal T(\M(D)) = 1 \mid s_0) = \Pr_\theta(\M(D) \in R \mid s_0) = \alpha,
    \]
    but this test has power
    \[
    \beta = \Pr_\theta( \mathcal T(\M(D)) = 1 \mid s_1) = \Pr_\theta(\M(D) \in R \mid s_1) > e^\eps \alpha .
    \]
\end{proofE}
Comparing Theorem~\ref{thm:pufferfish-hypothesis-tests} to Theorem~\ref{thm:dp-distingishing-neighbors}, we note strong parallels as well as some important differences. Both theorems provide interpretations for formal privacy guarantees (Pufferfish and DP, respectively) in terms of hypothesis tests, where sensitive level-$\alpha$ tests cannot have power exceeding $e^\eps \alpha$. The result for Pufferfish (Theorem~\ref{thm:pufferfish-hypothesis-tests}) is defined over random databases with flexible hypothesis tests, while the result for DP (Theorem~\ref{thm:dp-distingishing-neighbors}) is defined for fixed databases and specific hypothesis tests over complete, neighboring databases. We formally link the two concepts in Corollary~\ref{thm:dp-as-pufferfish}, which essentially states that DP is equivalent to a form of Pufferfish with a rigid definition of secrets but virtually no restrictions on the random process that generates a database.

\begin{corollaryE}
\label{thm:dp-as-pufferfish}
Let $\D$ denote a measurable space of databases, and let $\sim$ denote a neighboring relation between databases. Let $\mu(\D)$ denote the class of all measures over $\D$. Let $\sigma_D$ denote the predicate that the true database (in its entirety) is equal to $D$. Let $\Sset = \{\sigma_D \mid D \in \D \}$ denote these secrets, and let $\Spairs = \{ (\sigma_D, \sigma_{D'}) \mid D \sim D' \}$ denote the set of discriminative pairs (i.e., predicates involving neighboring databases).

Then $\M$ is $\eps$-DP($\D, \sim$) if and only if it is $\eps$-Pufferfish($\Sset, \Spairs, \mu(\D)$).
\end{corollaryE}

\begin{proofE}
    ($\implies$) Suppose the Pufferfish condition is satisfied. Fix $D_0 \sim D_1$, and let $\Pr_\theta$ be the probability measure on $\D$ that assigns probability $1/2$ to each of $D_0$ and $D_1$. Clearly we have $\Pr_\theta \in \mu(\D)$, and so for $D \sim \Pr_\theta$, Pufferfish implies:
    \begin{align*}
    \Pr(\M(D_0) \in E) &= \Pr_\theta(\M(D) \in E \mid D = D_0) \\
        &\leq e^\eps \, \Pr_\theta(\M(D) \in E \mid D = D_1) \\
        &= e^\eps \, \Pr(\M(D_1) \in E) .
    \end{align*}
    This is precisely the DP condition.

    ($\impliedby$) Suppose the Pufferfish condition is not satisfied. By Theorem~\ref{thm:pufferfish-hypothesis-tests}, there exist $D_0 \sim D_1$ from $\D$, $\Pr_\theta \in \mu(\D)$, $\alpha \in (0,1)$, and a test $\mathcal T(\Pr_\theta, \M, \M(D))$ such that for $D \sim \Pr_\theta$, $\mathcal T$ is level-$\alpha$ for $H_0: D = D_0$ vs. $H_1: D = D_1$, but $\mathcal T$ has power $\beta > e^\eps \alpha$.

    Define $\mathcal T_\theta(\cdot, \cdot) = \mathcal T(\Pr_\theta, \cdot, \cdot)$. Draw $D \sim \Pr_\theta$. By the above, $\mathcal T_\theta$ is level-$\alpha$ for $H_0: D = D_0$ vs. $H_1: D = D_1$ while having power $\beta > e^\eps \alpha$. Per Theorem~\ref{thm:dp-distingishing-neighbors}, $\M$ is not $\eps$-DP($\D, \sim$).
\end{proofE}
We highlight that Corollary~\ref{thm:dp-as-pufferfish} frames Pufferfish as a proper generalization of any pure DP variant. In particular, the instantiation of Pufferfish in Corollary~\ref{thm:dp-as-pufferfish} remains distribution-agnostic. This is a different presentation from in \citet[Section~6]{kifer2014pufferfish}, where various results equate specific DP settings with Pufferfish instantiations under the assumption of certain classes of data-generating distributions. While our result is general and unifying, it achieves this generality using the unnatural notion of complete, neighboring databases as the protected secrets. By contrast, the collection of results in \citet[Section~6]{kifer2014pufferfish} demonstrate that protection of individual tuples as secrets in classical DP requires the assumption of independent-tuple databases, in a manner quite similar to \citet{wasserman2010statistical}. (See Theorem~\ref{thm:wasserman} above.) Extending these results to the protection of individual network edges, as they note in a footnote, requires the assumption of independent edges.

\subsection{Edge Privacy in Pufferfish}
\label{sec:pufferfish-edge-dp}

It follows from Corollary~\ref{thm:dp-as-pufferfish} that edge DP corresponds to an instantiation of Pufferfish that protects against inference between complete, neighboring databases, but this provides no limits on the ability to infer specific edges in a network, as seen in Example~\ref{ex:insects}. In the style of \citet[Section~6]{kifer2014pufferfish}, we show that edge DP does protect against such inference in the presence of strong assumptions on the distribution of the network. In particular, the following theorem shows that edge DP protects against edge-level inference for networks with independent edges. (Note that for simplicity, the following theorem is stated for unattributed, undirected, binary graphs.)

\begin{theoremE}[][end]
\label{thm:edge-dp-as-pufferfish}
Let $\G_n$ be the set of binary graphs with vertices $V = [n]$, and let $\I_V$ denote a set of valid edges for $V$. Let $\Theta_{ind}$ be the set of all distributions on $\G_n$ whose edges occur independently, and let $\sigma_{ij}$ be the predicate that $\{i, j\} \in E$ for a given graph $G = (V, E)$. Denote the secrets
\[
\Sset = \{\sigma_{ij} : \{i,j\} \in \I_V \} \cup \{ \lnot \sigma_{ij} : \{i, j\} \in \I_V\}, \quad \Spairs = \{ (\sigma_{ij}, \lnot \sigma_{ij}) : \{i,j\} \in \I_V \}.
\]
Then $\M : \G_n \to \range(\M)$ is $\eps$--edge DP if and only if $\M$ is $\eps$--Pufferfish$(\Sset, \Spairs, \Theta_{ind})$.
\end{theoremE}

\begin{proofE}
The proof of this is similar to \citet[Theorem~6.1]{kifer2014pufferfish}.

$(\impliedby)$ Fix $g_0 \overset{e}{\sim} g_1 \in \G_n$. Then $g_0, g_1$ differ on one edge, say, $\{ i^*, j^* \} \not\in g_0$ and $\{ i^*, j^* \} \in g_1$. Let $\Pr_\theta$ be the distribution over $\G_n$ defined by
\[
\Pr_\theta(G = g) = \prod_{\{i,j\} \in \I_V} \pi_{ij}^{1_{\{i, j\} \in g}} (1 - \pi_{ij})^{1 - 1_{\{i, j\} \in g}} ,
\]
where
\[
\pi_{i^*j^*} = \frac 12, \quad \pi_{ij} = 1_{\{i, j\} \in g_0} \text{ for } \{i, j\} \in \I_V \setminus \{i^*, j^*\} .
\]
Note that $\Pr_\theta \in \Theta_{ind}$, and the support of $\Pr_\theta$ is precisely $\{ g_0, g_1 \}$. Specifically, $g_1$ is the only graph in the support of $\Pr_\theta$ satisfying $\sigma_{i^*j^*}$, and $g_0$ is the only graph satisfying $\lnot \sigma_{i^*j^*}$. So if $\M$ is $\eps$--Pufferfish$(\Sset, \Spairs, \Theta_{ind})$ we have:
\[
\Pr_\theta(\M(G) = \omega \mid \sigma_{i^* j^*}) \leq e^\eps \Pr_\theta(\M(G) = \omega \mid \lnot \sigma_{i^* j^*}),
\]
or equivalently:
\[
\Pr(\M(g_1) = \omega) \leq e^\eps \Pr(\M(g_0) = \omega).
\]
A similar argument shows $\Pr(\M(g_0) = \omega) \leq e^\eps \Pr(\M(g_1) = \omega)$, and so we have $\eps$--edge DP, since $g_0$ and $g_1$ are arbitrary neighboring databases in $\G_n$.

$(\implies)$ Let $\{ i^*, j^* \} \in \I_V$, $\Pr_\theta \in \Theta_{ind}$.
%We will distinguish between $G$ random and $g$ fixed with upper- and lowercase notation.
%For a fixed graph $g$, let $g \lor 1$ denote $g$ taken to have $A(g)_{i^*j^*} = 1$ and $g \land 0$ denote $g$ with $A(g)_{i^*j^*} = 0$.
Assuming $\M$ is $\eps$--edge DP, we have:
\[
\begin{aligned}
\Pr_\theta(\M(G) = \omega \mid \sigma_{i^*j^*}) 
    &= \sum_{g \in \G_n} \Pr(\M(g) = \omega) \Pr_\theta(G = g \mid \sigma_{i^*j^*}) \\
    &= \sum_{g \in \G_n} \Pr(\M(g \cup \{ i^*, j^* \}) = \omega) \Pr_\theta(G = g \mid \sigma_{i^*j^*}) \quad\text{($g = g \cup \{i^*, j^* \}$ under $\sigma_{i^*j^*}$)} \\
    &\leq e^\eps \sum_{g \in \G_n} \Pr(\M(g \setminus \{i^*, j^* \}) = \omega) \Pr_\theta(G = g \mid \sigma_{i^*j^*})  \quad\text{(by edge DP)} \\
    \text{\textcircled{a}} &= e^\eps \sum_{g \in \G_n} \Pr(\M(g \setminus \{i^*, j^* \}) = \omega) \Pr_\theta(G = g \mid \lnot \sigma_{i^*j^*}) \\
    &= e^\eps \sum_{g \in \G_n} \Pr(\M(g) = \omega) \Pr_\theta(G = g \mid \lnot \sigma_{i^*j^*}) \quad\text{($g = g \setminus \{i^*, j^* \}$ under $\lnot \sigma_{i^*j^*}$)} \\
    &= e^\eps \Pr_\theta(\M(G) = \omega \mid \lnot \sigma_{i^* j^*}) .
\end{aligned}
\]
Equality \textcircled{a} comes from the fact that for any $\Pr_\theta \in \Theta_{ind}$, the set
\[
\{ g \in \G_n : \Pr_\theta(G = g \mid \sigma_{i^*j^*}) > 0 \}
\]
has a one-to-one correspondence with the set
\[
\{ g' \in \G_n : \Pr_\theta(G = g' \mid \lnot \sigma_{i^*j^*}) > 0 \}
\]
such that $g \mapsto g \setminus \{i^*, j^* \}$, and
\[
\Pr_\theta(G = g \mid \sigma_{i^* j^*}) = \prod_{\{i,j\} \in \I_V \setminus \{i^*, j^*\}} f_{ij}(1_{\{i, j\} \in g}) = \Pr_\theta(G = g' \mid \lnot \sigma_{i^* j^*}) ,
\]
where each $f_{ij} : \{0, 1\} \to \R_{\geq 0}$ is some function whose particular form is not of concern.

By a similar argument, we have $\Pr_\theta(\M(G) = \omega \mid \lnot \sigma_{i^*j^*}) \leq e^\eps \Pr_\theta(\M(G) = \omega \mid \sigma_{i^* j^*})$. Since $i^*, j^*$ are arbitrary, we meet the $\eps$-Pufferfish constraints for $(\Sset, \Spairs, \Theta_{ind})$.
\end{proofE}
Theorem~\ref{thm:edge-dp-as-pufferfish} states the result that we've come to expect from the preceding discussion: edge DP \textit{does} limit edge disclosure risk on networks with independent edges. In other words, the misinterpretations of edge DP highlighted in Section~\ref{sec:interpretations}---i.e., broad claims that edge DP limit edge disclosure risk---do hold in highly specialized cases. In sum, we can say that for all networks, an adversary's ability to distinguish between two complete, neighboring databases is hindered by edge DP (Corollary~\ref{thm:dp-as-pufferfish}), while for independent-edge networks, we can additionally make the more desirable claim that an adversary's ability to infer specific edges is equally hindered by edge DP. An outstanding question is whether for other sorts of networks we can make a weaker claim on the adversary's ability to infer specific edges. Lemma~\ref{thm:edge-dp-plus-alpha-pufferfish} provides one way to approach this.

\begin{lemmaE}[][end]
\label{thm:edge-dp-plus-alpha-pufferfish}
In the setting of Theorem~\ref{thm:edge-dp-as-pufferfish}, let $\Theta$ denote any set of distributions over $\G_n$.
Assume further that there exists $\alpha \in [0, \infty)$ such that for all $\Pr_\theta \in \Theta$, $g \in \G_n$, and $\{i, j\} \in \I_{V}$:
\begin{align}
    \label{eq:alpha-pufferfish}
    \Pr_\theta(G = g \cup \{i, j\} \mid \{i, j\} \in E) &\leq e^\alpha \Pr_\theta(G = g \setminus \{i,j\} \mid \{i, j\} \not\in E) \\
    \Pr_\theta(G = g \setminus \{i, j\} \mid \{i, j\} \not\in E) &\leq e^\alpha \Pr_\theta(G = g \cup \{i,j\} \mid \{i, j\} \in E) . \nonumber
\end{align}
Then if $\M$ is $\eps$--edge DP, it is also $(\eps + \alpha)$--Pufferfish$(\Sset, \Spairs, \Theta)$.
\end{lemmaE}

\begin{proofE}
To satisfy $(\eps + \alpha)$--Pufferfish$(\Theta, \Sset, \Spairs)$, we need to satisfy
\begin{align*}
\Pr_\theta(\M(G) = \omega \mid \{ i, j \} \in E ) &\leq e^{\eps + \alpha} \Pr_\theta(\M(G) = \omega \mid \{ i, j \} \not\in E) \\
\Pr_\theta(\M(G) = \omega \mid \{ i, j \} \not\in E) &\leq e^{\eps + \alpha} \Pr_\theta(\M(G) = \omega \mid \{ i, j \} \in E)
\end{align*}
for all $\Pr_\theta \in \Theta, \{i, j\} \in \I_V, \omega \in \range(\M)$ such that $\Pr_\theta(\{i, j\} \in E) \in (0, 1)$. Fix $\Pr_\theta, i, j$, and $\omega$. We begin by rewriting
\begin{align*}
    \Pr_\theta(\M(G) = \omega \mid \{ i, j \} \in E ) &= \sum_{g \in \G_n} \Pr(\M(g) = \omega) \Pr_\theta(G = g \mid \{i, j\} \in E) \\
        &= \frac 12 \sum_{g \in \G_n} \Pr(\M(g \cup \{i, j\}) = \omega) \Pr_\theta(G = g \cup \{i, j\} \mid \{i, j\} \in E),
\end{align*}
where the last line follows from the fact that replacing $g$ with $g \cup \{i, j\}$ results in double-counting, as $\Pr_\theta(G = g \mid \{i, j\} \in E) = 0$ for all $g$ not containing the edge $\{i, j\}$. We similarly rewrite:
\[
    \Pr_\theta(\M(G) = \omega \mid \{ i, j \} \not\in E ) = \frac 12 \sum_{g \in \G_n} \Pr(\M(g \setminus \{i, j\}) = \omega) \Pr_\theta(G = g \setminus \{i, j\} \mid \{i, j\} \not\in E) .
\]

For any fixed $g$, since $\M$ is $\eps$--edge DP, we know that
\begin{align*}
\Pr(\M(g \cup \{i, j\}) = \omega) &\leq e^\eps \Pr(\M(g \setminus \{i, j\}) = \omega) , \text{ and} \\
\Pr(\M(g \setminus \{i, j\}) = \omega) &\leq e^\eps \Pr(\M(g \cup \{i, j\}) = \omega) ,
\end{align*}
since $g \cup \{i, j\}$ neighbors $g \setminus \{i, j\}$, while by assumption we also have:
\begin{align*}
    \Pr_\theta(G = g \cup \{i, j\} \mid \{i, j\} \in E) &\leq e^\alpha \Pr_\theta(G = g \setminus \{i,j\} \mid \{i, j\} \not\in E) , \text{ and} \\
    \Pr_\theta(G = g \setminus \{i, j\} \mid \{i, j\} \not\in E) &\leq e^\alpha \Pr_\theta(G = g \cup \{i,j\} \mid \{i, j\} \in E) .
\end{align*}
Consequently we may conclude,
\begin{align*}
    &\Pr_\theta(\M(G) = \omega \mid \{ i, j \} \in E ) \\
    &\quad = \frac 12 \sum_{g \in \G_n} \Pr(\M(g \cup \{i, j\}) = \omega) \Pr_\theta(G = g \cup \{i, j\} \mid \{i, j\} \in E) \\
    &\quad \leq \frac 12 \sum_{g \in \G_n} \Big( e^\eps \Pr(\M(g \setminus \{i, j\}) = \omega) \Big) \, \Big( e^\alpha \Pr_\theta(G = g \setminus \{i, j\} \mid \{i, j\} \not\in E) \Big) \\
    &\quad= e^{\eps + \alpha} \Pr_\theta(\M(G) = \omega \mid \{ i, j \} \not\in E ) ,
\end{align*}
and similarly:
\[
    \Pr_\theta(\M(G) = \omega \mid \{ i, j \} \not\in E ) \leq e^{\eps + \alpha} \Pr_\theta(\M(G) = \omega \mid \{ i, j \} \in E ) .
\]
\end{proofE}
Intuitively, Lemma~\ref{thm:edge-dp-plus-alpha-pufferfish} requires that if we condition on the presence or absence of a single edge, the distribution of the remaining edges should not significantly change. This is precisely the case for independent-edge graphs, where we obtain $\alpha = 0$ to recover the results of Theorem~\ref{thm:edge-dp-as-pufferfish}. At the other end of the spectrum lies our insect example (Example~\ref{ex:insects}). In this example, the presence or absence of an edge between the two queens determines the overall density of the remainder of the network. As $n$ grows, the inequalities in Eq.~\ref{eq:alpha-pufferfish} hold only when $\alpha$ tends to $\infty$. Edge DP here offers no guarantee for protecting disclosure of the queens' edge.

\subsection{Application to Exponential Random Graph Models}

For many families of distributions, the conditions of Lemma~\ref{thm:edge-dp-plus-alpha-pufferfish} are difficult to verify. One convenient class of distributions, however, is the family of \textit{exponential random graph models}, or \textit{ERGMs} \citep{wasserman1996logit,snijders2006new}, for which we show how to apply Lemma~\ref{thm:edge-dp-plus-alpha-pufferfish}.

An exponential random graph model $\Pr_\theta \in \Theta_{ergm}$ is defined as an exponential-family distribution over the set of possible graphs $g \in \G_n$,
\[
\Pr_\theta(G = g) = \frac{\exp(\beta^T u(g))}{\kappa(\beta)} ,
\]
where $\beta \in \R^d$ is a $d$-dimensional vector of parameters, $u(g) \in \R^d$ is a vector of sufficient network statistics, and $\kappa(\beta)$ is a normalizing constant. Generally speaking, $u(g)$ may flexibly accommodate a wide variety of statistics, including the total number of edges (to capture network density), the number of mutual edges (to capture reciprocity in the case of directed graphs), the number of edges among nodes sharing a given covariate value (to capture homophily in the case of attributed graphs), various structural features such as the number of $k$-stars, triangles, etc., among many other examples. We restrict our attention, however, to statistics arising from the network edges only, as a matter of both simplicity and practicality, since we are concerned with edge privacy.%
\footnote{In the case of attributed networks, edge DP alone offers no privacy for node attributes, as the definition of edge neighbors ignores node attributes. Some practitioners use edge DP on these graphs anyway, warning that node attributes are not protected and are assumed to be publicly known \citep[e.g.,][]{karwa2017sharing}. For more information, see the discussion in Section~\ref{sec:conclusion}.}

It is well-known that ERGMs admit a convenient logistic form, where the odds that any edge occurs, conditioned on the remainder of the graph, can be written
\begin{align}
\label{eq:change-statistic}
\frac{\Pr_\theta( \{ i, j \} \in E \mid G \cup \{i, j\} = g \cup \{i, j \})}{\Pr_\theta( \{ i, j \} \not\in E \mid G \setminus \{i, j\} = g \setminus \{i, j \})}
    &= \frac{\exp(\beta^T u(g \cup \{i, j\})}{\exp(\beta^T u(g \setminus \{i, j\})} \\
    &= \exp(\beta^T \Delta(g,i,j)) , \nonumber
\end{align}
where $\Delta(g,i,j) = u(g \cup \{i, j\}) - u(g \setminus \{i, j\})$ is the \textit{change statistic} for a graph $g \in \G_n$ \citep{hunter2006inference,hunter2012computational}. This nearly gives us the constant $\alpha$ from Lemma~\ref{thm:edge-dp-plus-alpha-pufferfish}. The final step is to apply Bayes' theorem to reverse the conditional probabilities. This is reflected in Corollary~\ref{thm:ergm}.

\begin{corollaryE}
\label{thm:ergm}
Let $\Theta \subseteq \Theta_{ergm}$ denote a family of ERGMs on $n$ nodes using edge-based statistics. Let $\beta_\theta$ denote the ERGM parameters for a given model $\Pr_\theta \in \Theta$, and let $\Delta_\theta(g,i,j)$ denote a change statistic for a given model on a given edge and a given realized graph. Then any privacy mechanism $\M$ that is $\eps$--edge DP is also $(\eps + \alpha)$--Pufferfish$(\Theta, \Sset, \Spairs)$ for the edge secrets $\Sset, \Spairs$ from Theorem~\ref{thm:edge-dp-as-pufferfish} and
\[
\alpha = \sup_{\theta, g, i, j} \left| \, \beta_\theta^T \Delta_\theta(g, i, j)) - \log \left( \frac{\Pr_\theta( \{i, j \} \in E)}{\Pr_\theta( \{i, j\} \not\in E)} \right) \right| \leq \sup_{\theta, g, i, j} \left| 2 \beta_\theta^T \Delta_\theta(g, i, j) \right| .
\]
\end{corollaryE}

\begin{proofE}
The derivation of $\alpha$ follows from an application of Bayes' theorem to the definition of the change statistic in Eq.~\ref{eq:change-statistic}. In particular, finding $\alpha$ to satisfy Lemma~\ref{thm:edge-dp-plus-alpha-pufferfish} requires us to bound the quantity
\[
\left| \log \left( \frac{\Pr_\theta( G = g \cup \{i, j \} \mid \{i, j\} \in E)}{\Pr_\theta( G = g \setminus \{i, j \} \mid \{i, j\} \not\in E)} \right) \right|.
\]
To simplify notation, we write the event $(G \cup \{i, j\} = g \cup \{i, j\}) = (G \setminus \{i, j\} = g \setminus \{i, j\})$ as $G_{ij}^c = g_{ij}^c$. By Bayes' theorem, we have that:
\begin{align*}
\Pr_\theta( G = g \cup \{i, j \} \mid \{i, j\} \in E)
    &= \Pr_\theta( G \cup \{i, j\} = g \cup \{i, j \} \mid \{i, j\} \in E) \\
    &= \Pr_\theta(G_{ij}^c = g_{ij}^c \mid \{i, j\} \in E) \\
    &= \frac{\Pr_\theta(\{i, j\} \in E \mid G_{ij}^c = g_{ij}^c) \Pr_\theta(G_{ij}^c = g_{ij}^c)}{\Pr_\theta(\{i, j\} \in E)},
\end{align*}
and similarly
\[
\Pr_\theta( G = g \setminus \{i, j \} \mid \{i, j\} \not\in E) = \frac{\Pr_\theta(\{i, j\} \not\in E \mid G_{ij}^c = g_{ij}^c) \Pr_\theta(G_{ij}^c = g_{ij}^c)}{\Pr_\theta(\{i, j\} \not\in E)} .
\]
Dividing these two quantities, we obtain:
\[
\frac{\Pr_\theta( G = g \cup \{i, j \} \mid \{i, j\} \in E)}{\Pr_\theta( G = g \setminus \{i, j \} \mid \{i, j\} \not\in E)} = \left. \underbrace{\frac{\Pr_\theta(\{i, j\} \in E \mid G_{ij}^c = g_{ij}^c)}{\Pr_\theta(\{i, j\} \not\in E \mid G_{ij}^c = g_{ij}^c)}}_{\textcircled{a}} \middle/ \underbrace{\frac{\Pr_\theta(\{i, j\} \in E)}{\Pr_\theta(\{i, j\} \not\in E)}}_{\textcircled{b}} \right. .
\]
By Eq.~\ref{eq:change-statistic}, $\textcircled{a}$ is precisely $\exp(\beta_\theta^T \Delta_\theta(g, i, j))$. Thus:
\[
\left| \log \left( \frac{\Pr_\theta(\{i, j\} \in E \mid G_{ij}^c = g_{ij}^c)}{\Pr_\theta(\{i, j\} \not\in E \mid G_{ij}^c = g_{ij}^c)} \middle/ \frac{\Pr_\theta(\{i, j\} \in E)}{\Pr_\theta(\{i, j\} \not\in E)} \right) \right| = \left| \beta_\theta^T \Delta_\theta(g, i, j) - \log\left( \frac{\Pr_\theta(\{i, j\} \in E)}{\Pr_\theta(\{i, j\} \not\in E)} \right) \right| .
\]
Taking the supremum over all choices of $\theta, g, i, j$ yields $\alpha$.

We note that $\log(\textcircled{b})$ is also bounded in absolute value by a similar quantity, yielding the upper bound on $\alpha$. Let $\zeta = \sup_{\theta, g, i, j} \left| \beta_\theta^T \Delta_\theta(g, i, j) \right|$. Then:
\begin{align*}
\frac{\Pr_\theta( \{i, j\} \in E)}{\Pr_\theta( \{i, j\} \not\in E)}
    &= \frac{\sum_{g \in \G_n} \Pr_\theta( \{i,j\} \in E \mid G_{ij}^c = g_{ij}^c) \Pr_\theta(G_{ij}^c = g_{ij}^c)}{\sum_{g \in \G_n} \Pr_\theta \{i,j\} \not\in E \mid G_{ij}^c = g_{ij}^c) \Pr_\theta(G_{ij}^c = g_{ij}^c)}.
\end{align*}
It follows from Eq.~\ref{eq:change-statistic} that for each $g \in \G_n$:
\[
\Pr_\theta( \{i,j\} \in E \mid G_{ij}^c = g_{ij}^c) \leq e^{\zeta} \Pr_\theta( \{i,j\} \not\in E \mid G_{ij}^c = g_{ij}^c),
\]
and similarly
\[
\Pr_\theta( \{i,j\} \in E \mid G_{ij}^c = g_{ij}^c) \geq e^{-\zeta} \Pr_\theta( \{i,j\} \not\in E \mid G_{ij}^c = g_{ij}^c).
\]
Plugging these bounds into the previous equation yields:
\[
e^{-\zeta} \leq \frac{\Pr_\theta( \{i, j\} \in E)}{\Pr_\theta( \{i, j\} \not\in E)} \leq e^\zeta \quad \iff \quad \left| \log \left( \frac{\Pr_\theta( \{i, j \} \in E)}{\Pr_\theta( \{i, j\} \not\in E)}  \right) \right| \leq \zeta.
\]
Putting it all together, we have that:
\[
\alpha = \sup_{\theta, g, i, j} | \log(\textcircled{a}) - \log(\textcircled{b}) | \leq 2 \zeta .
\]
\end{proofE}
Corollary~\ref{thm:ergm} provides both an explicit form of $\alpha$ as well as a simplified upper bound on $\alpha$. Generally speaking, when $\beta_\theta^T \Delta_\theta(g, i, j)$ is bounded, an $\eps$--edge DP mechanism does provide some protection against inferring individual edges, although possibly weaker than implied by $\eps$ alone. On the other hand, when this quantity is unbounded, an $\eps$--DP mechanism may fail to protect against such inference.

\begin{remark}
    Consistent with Theorem~\ref{thm:edge-dp-as-pufferfish}, Corollary~\ref{thm:ergm} yields a value of $\alpha = 0$ for independent-edge ERGMs. This is because for any $\Pr_\theta \in \Theta_{ergm} \cap \Theta_{ind}$, the marginal log-odds of a given edge occurring are identical to the log-odds conditioned on the remainder of the graph, as given in Eq.~\ref{eq:change-statistic}.
\end{remark}

\section{Discussion and Open Questions}
\label{sec:conclusion}

Nearly twenty years after its debut in the literature, researchers are still trying to find clear ways to communicate the meaning of differential privacy to both data users and data participants \citep{cummings2021need,nanayakkara2023chances,cummings2023challenges}. Quoting \citet{nanayakkara2023chances}:
\begin{quote}
    [I]t is vital that organizations deploying DP effectively communicate the privacy implications of implementation details that govern the strength of systems' privacy protections. Without such transparency, organizations risk engaging in ``privacy theater,'' which may result in people falsely believing they are well-protected.
\end{quote}
We contend that this risk is perhaps especially important in the context of network data, where dyadic relationships make causal interpretations of privacy difficult, and dependence structures pose major challenges for inferential interpretations. Given the abundance of inaccurate interpretations of network DP in the research literature, it stands to reason that effective communication of network DP guarantees to end-users is in need of special attention.

We have shown in particular how edge DP has been misinterpreted in the literature, leaning on inferential interpretations that lack sufficient qualification. Moreover, we have demonstrated some theoretical conditions under which the gap can be bridged between these misinterpretations and the guarantees actually implied by edge DP. A major challenge remains to bridge the gap between these theoretical developments and effective communication to end-users, but we believe clarifying the theoretical understanding among researchers is an important first step.

The primary theoretical focus of this work has been the application of pure edge DP to unattributed, binary-valued graphs. While many of the results stated generalize to directed or undirected graphs, with or without self-loops, some important theoretical gaps remain. We enumerate several open questions here.

\textbf{Attributed Graphs and Node Privacy.} In the case of attributed graphs (i.e., networks with node covariates) edge DP alone offers no privacy for node attributes, since the privacy definition does not specify how a mechanism must treat differences in node-level data. Treatment of attributed graphs in practice varies, sometimes proceeding with edge DP anyway and other times choosing alternative privacy definitions. Edge DP's primary competitor is node DP, which, informally, defines two networks as neighbors if they differ on a single node or any of its incident edges.\footnote{As in edge DP, there are variations on how exactly this is defined, with the original definition of \citet{hay2009accurate} allowing for the addition or removal of a node, and others allowing only a substitution of one node with another.} This provides a stronger guarantee that would include protection of node attribute information. Node DP is arguably better suited than edge DP for causal-flavored interpretations, since two hypothetical neighboring databases may correspond to, e.g., a network in which a given node has consented to share his node and edge data vs. one in which the node has not consented to such sharing. But the strength of node DP's privacy guarantee often limits its use in practice \citep{kasiviswanathan2013analyzing}.

It should be noted that edge DP and node DP are neither mutually exclusive, nor are they the only DP definitions applied to networks. The quest to balance privacy and utility in network analysis has led to a number of other approaches, including hybrid definitions such as \textit{edge-adjacent DP} \citep{blocki2013differentially} and the composition of separate edge DP and node DP mechanisms for edge vs. node attribute data \citep[e.g.,][]{ji2019differentially}.%
\footnote{It can be shown that the composition of separate edge and node DP mechanisms in this manner yields an edge-adjacent DP guarantee.}
These hybrid definitions raise further questions for interpretation of the privacy guarantee. For example, in edge-adjacent DP, two networks are defined as neighbors if they differ on either a single node's attributes or a single edge in the network.%
\footnote{Once again, the literature contains inconsistent usage of this privacy definition. In \citet{blocki2013differentially}, only one of these changes is permitted between neighboring networks, while in \citet{jorgensen2016publishing,chen2019publishing} both changes are allowed simultaneously.}
This definition of neighbors leads to another situation that is difficult to interpret causally: what privacy situation leads to a network in which either a single node's attributes or a single edge has been manipulated? Meanwhile, dependence between node attributes and edge formation leads to additional challenges in providing inferential-privacy guarantees.

\textbf{Approximate DP.} The definitions of DP used in this work cover only $\eps$-DP (i.e., \textit{pure DP}), in contrast with approximate forms of differential privacy like $(\eps, \delta)$--DP, zero-concentrated DP, Renyi DP, and $f$-DP. In contrast with pure DP, these approximate forms of DP may provide a weakened privacy guarantee, allowing for the development of privacy mechanisms that retain additional data utility. Generally speaking, interpretation of DP approximations is harder than interpretation of pure DP. In fact, the $\delta$ parameter of the popular $(\eps, \delta)$--DP formulation is frequently misinterpreted, and obtaining a correct interpretation is a fairly technical endeavor. (See discussion and citations in \citealp{kifer2022bayesian}, Section~5.1.) Generalizing inferential-privacy results for approximate forms of network DP remains an open question. Of the various approximations in circulation, $f$-DP \citep{dong2022gaussian} is a natural starting point, since it is parameterized in terms of a trade-off function $f$ between the Type I and Type II errors of the hypothesis test for neighboring databases $H_0: D = D_0$ vs. $H_1: D = D_1$.

\textbf{Data Collection and Local Differential Privacy.} A particularly challenging aspect of network privacy is the question of who owns the data and how it is collected. The privacy definitions we have focused on in this work all fall under the umbrella of \textit{central DP}, where one trusted party is assumed to have access to the full network data, and that party applies the privacy mechanism to the data to produce a sanitized release. By contrast, \textit{local DP} \citep{duchi2013local} describes the situation in which data is collected from distributed parties who each perturb their data using a privacy mechanism. These individual local DP ``views'' may then be analyzed separately or in aggregate. The application of local DP to network data is more nuanced than in the traditional setting. For example, in an undirected graph, what party is responsible for reporting a given edge? This question is at the heart of the distinction between the closely-related definitions of edge-local DP \citep{qin2017generating} and relationship DP \citep{imola2021locally}.
Proper interpretations of DP guarantees for network data should account for these nuances in data collection and privacy definition.

\appendix

\section{Proofs}

In the proofs to follow, we denote a hypothesis test for $H_0$ vs. $H_1$ as a (possibly) randomized algorithm $\mathcal T : \range(\M) \to \{0, 1\}$, where we reject $H_0$ iff $\mathcal T(\M(D)) = 1$.

\printProofs

\printbibliography

\end{document}